\documentclass[reqno]{amsart}
\usepackage{epsfig}
\usepackage{amsmath}
\usepackage{amsthm}
\usepackage{amssymb}
\usepackage{hodge}
\usepackage{graphicx}
\usepackage{graphics}
\usepackage{psfrag}
\usepackage{xypic}

\newcommand{\be}{\begin{equation}}
\newcommand{\ee}{\end{equation}}
\numberwithin{equation}{section}

\newcommand{\CM}{{\mathcal M}}
\newcommand{\CKM}{{\mathcal K}{\mathcal M}}
\newcommand{\re}{{\rm e}}
\newcommand{\CF}{{\mathcal F}}

\begin{document}

\title{On the mathematics and physics of high genus invariants of $\orb$}

\author{Vincent Bouchard}
\address{Harvard University, Jefferson Physical Laboratory, 
17 Oxford St., Cambridge, MA 02138, USA}
\email{bouchard@physics.harvard.edu}
\author{Renzo Cavalieri}
\address{University of Michigan, Department of Mathematics,
2074 East Hall, 530 Church Street, Ann Arbor, MI 48109-1043
}
\email{crenzo@umich.edu}


\begin{abstract}
This paper wishes to foster communication between  mathematicians and  physicists working in mirror symmetry and orbifold \GW\ theory. We provide a reader friendly  review of the physics computation in \cite{ABK} that  predicts  \GW\ invariants of $\orb$ in arbitrary genus, and of the mathematical framework for expressing these invariants as Hodge integrals. 
Using geometric properties of the Hodge classes, we compute the unpointed invariants for $g=2,3$, thus providing the first high genus mathematical check of the physics predictions.

\end{abstract}
\maketitle

\section*{Introduction}
\subsection*{Scope and results}

All too often mathematicians and physicists are compared to a couple in a disfunctional marriage: sharing a household but unable to communicate properly. This paper attempts to contradict this stereotype, by exploring the orbifold \GW\ theory of $\orb$.

On the one hand, we distill for a mathematical audience, in sections 3 and 4, the physics calculation of \cite{ABK}, which provides a prediction for the unmarked and marked Gromov-Witten invariants of $\orb$ at any genus. This calculation is close in spirit to the original calculation of the number of rational curves in the quintic threefold by Candelas et al \cite{Candelas}, relying on mirror symmetry and topological string theory. 

On the other hand, $\orb$ invariants can be defined mathematically, and interpreted in terms of \zhh: top intersections of characteristic classes of some natural vector bundle on  moduli spaces of covers of curves. We present this point of view in section 1, trying to cater especially to the physicist reader. \zhh\ are new mathematical creatures, and their systematic exploration is on the second author's research agenda. In low genus, some ad hoc considerations lead to the following original result, which is proved in section 2.

\begin{thm}
The unpointed invariants of $\orb$ are mathematically computed for $g=2$ and $3$, and agree with the predictions of \cite{ABK}.
\end{thm}

This result provides an interesting validity check of the high genus predictions of \cite{ABK}, since so far only the genus $0$ predictions had been proved mathematically, computed in three independent ways by Coates, Corti, Iritani and Tseng in \cite{CCIT}, Bayer and Cadman in \cite{bc:c3z3} and by Cadman and Cavalieri \cite{cc:c3z3}.

\subsection*{History and connections}

The orbifold $\orb$ has recently been an exciting object of study both for mathematicians and physicists. In mirror symmetry, it represents a special point in the stringy K\"ahler moduli space of its crepant resolution, local $\Proj$. This point of view has been used in various ways in the past to study string physics on $\orb$ --- see for instance \cite{OFS,DG} for D-brane aspects. However, only recently was it used to relate the (orbifold) \GW\ theory of $\orb$ to the \GW\ theory of local $\Proj$ \cite{ABK}. Mathematically, this is an incarnation of the \textit{McKay philosophy}, stating that the $G$-equivariant geometry of a space $X$ should equal the geometry of a crepant resolution of the quotient $X/G$. 
Precise statements about this equivalence in \GW\ theory have been formulated by Ruan (\cite{r:crc}), Bryan and Graber under some technical assumption on the target orbifold (\cite{bg:crc}) and Coates, Corti, Iritani and Tseng (\cite{ccit:crcas}). 
These conjectures have been verified in genus $0$ for several examples(\cite{bg:crc}, \cite{ccit:wc}, \cite{ccit:crcas}, \cite{CCIT} \cite{g:crc}, \cite{bg:crchhi}, \cite{bg:root}). To the best of our knowledge, no examples for higher genus have been worked out yet.

\subsection*{Acknowledgments}

We owe  special thanks to M. Aganagic: it is through a three-way  exchange of ideas while the authors were visiting Berkeley that this project got started.
We would also like to thank J. Bryan, T. Coates, A. Klemm and R. Vakil for interesting discussions. Preprint number: HUTP-07/A0005.

\section{The mathematics}

We first review some aspects of the mathematics of orbifold Gromov-Witten invariants, and then focus on the main character $\orb$ of this note.

\subsection{Orbifold Gromov-Witten invariants}

Let $\mathfrak{X}$ be an orbifold, or, if you prefer, a Deligne Mumford stack.
The study of \GW\ invariants of orbifolds is developed by Chen and Ruan in \cite{cr:ogwt} and \cite{cr:nctoo}. The algebraic point of view is established in \cite{Abramovich_Graber_Vistoli}. In order to obtain a good mathematical theory (i.e. a compact and reasonably well behaved moduli space, equipped with a virtual fundamental class) they introduce the following two modifications to the ordinary \GW\ set-up:
\begin{description}
	\item[twisted stable maps] the source curves must be allowed to become ``stacky''. Informally, a twisted stable curve is ``almost'' a curve: it has a finite set of twisted points, where it locally looks like $[\bC/\bZ_n]$, the (stack) quotient of $\bC$ by the action of a cyclic group. Ordinary stable maps are replaced by (representable) morphisms from twisted stable curves. 
	\item[orbifold cohomology insertions] ordinary \GW\ invariants have insertions that take value in the cohomology of the target space.  Here, one needs to enlarge cohomology to the Chen-Ruan orbifold cohomology ring, including classes that contain a combination of geometric and representation theoretic data, keeping track of the automorphisms that the cohomology classes might have. Formally, this is defined to be the cohomology of a related orbifold $I\mathfrak{X}$, called the inertia orbifold.  
\end{description}
With these two modifications in place, the moduli space $\overline{M}_{g,n}(\mathfrak{X},\beta)$ is a proper Deligne Mumford stack of expected dimension
\be
(1-g)(\mbox{dim}\ \mathfrak{X}-3) - K_\mathfrak{X}\cdot \beta +n ,
\ee
and just about any desirable (and undesirable) feature of ordinary \GW\ theory carries over to the orbifold setting.

\subsection{Twisted stable maps to $\mathcal{B}\mathbb{Z}_3$}

Consider the orbifold $\mathfrak{X}= \bz$, which can be thought of as the classifying space for principal $\bZ_3$ bundles, or as the global quotient $[pt/\bZ_3]$ of a point by the trivial action of the group $\bZ_3$.

\begin{figure}
\begin{center}
\includegraphics[width=5in]{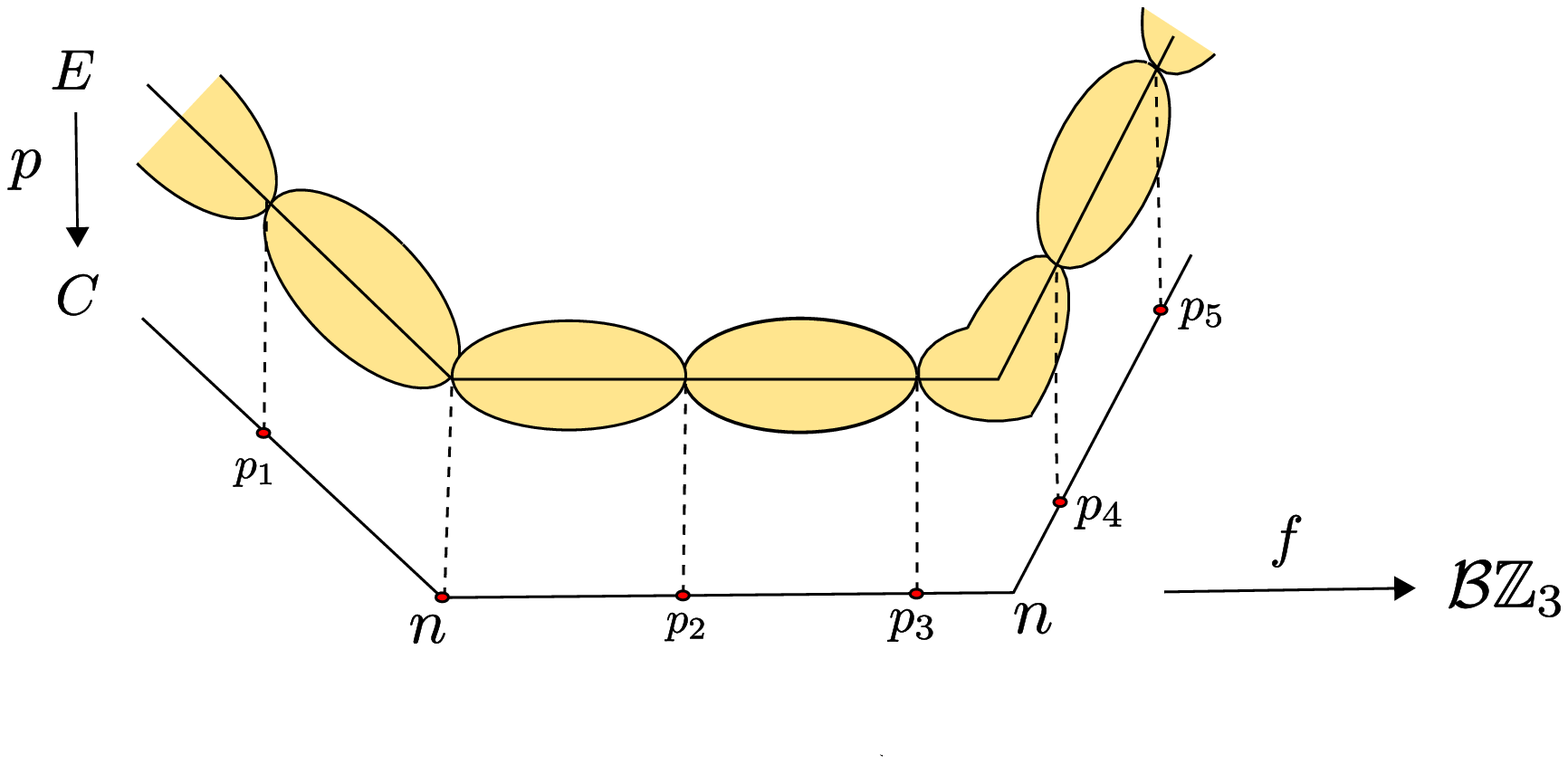}
\end{center}
\caption{A map $\mathcal{C} \rightarrow \bz$ corresponds to an admissible cover of $\mathcal{C}$. Notice that not all marked points must be twisted.}
\label{fig:maptobz3}
\end{figure}

In \cite{Abramovich_Corti_Vistoli}, Abramovich, Corti and Vistoli show that the stack $\sm{n}{\bz}{0}$ is the (normalization of the) moduli space of admissible $\mathbb{Z}_3$-covers of genus $g$ curves. This stack parameterizes degree $3$ covers $p:E\rightarrow C$ such that:
\begin{itemize}
	\item $C$ is a stable $(n)$-marked genus $g$ curve (the coarse moduli space of the twisted curve $\mathcal{C}$);
	\item $E$ is a nodal curve; nodes of $E$ ``correspond to''\footnote{The preimages of nodal (resp. smooth) points of $C$ are nodal (resp. smooth) points of $E$.} nodes of $C$;
	\item $E$ is endowed with a $\mathbb{Z}_3$ action;
	\item $p$ is the quotient map with respect to the action;
	\item $p$ is ramified only over the marked points of $C$, and possibly over the nodes;
	\item when $p$ is ramified over a node, denote $x_1$ and $x_2$ the shadows of the node in the normalization $\tilde{E}$. The $\mathbb{Z}_3$-representations induced on $T_{x_1}$ and $T_{x_2}$ are dual to each other. 
\end{itemize}
This description is illustrated in Figure \ref{fig:maptobz3}.


We turn our attention now to the case with no marks. A general point in the moduli space $\smu{\bz}{0}$ represents an \'{e}tale $\bZ_3$-cover of a smooth genus $g$ curve $C$, equivalent to the data of the curve $C$ and a monodromy representation, i.e. a group homomorphism
\be
\varphi: \pi_1(C) \rightarrow \bZ_3.
\ee

The canonical forgetful morphism
\be
\smu{\bz}{0} \longrightarrow \overline{M}_g
\ee
is finite of degree\footnote{The factor of $1/3$ comes from the fact that every cover has a degree $3$ non-trivial automorphism, given by the action of a generator of $\bZ_3$.} $3^{2g}/3$, but not \'{e}tale: it ramifies over the boundary of $\overline{M}_g$.

It is important to observe that the moduli space $\smu{\bz}{0}$ consists of two connected components:
\begin{description}
	\item[$\smu{\bz}{0}^{disc}$] parameterizes disconnected covers: three copies of $C$ mapping down to $C$ via  identity maps. These covers correspond to the trivial monodromy. This component is in fact essentially a copy of $\overline{M}_g$: the only difference is that the covers have a degree $3$ nontrivial automorphism. Therefore the forgetful map restricted to this component has degree $1/3$. 
	\item[$\smu{\bz}{0}^{conn}$]	parameterizes connected covers, corresponding to nontrivial monodromy representations.
	\end{description}

\subsection{Hodge bundles}
The Hodge bundle $\mathbb{E}^h$ is a rank $h$ vector bundle on $\overline{M}_h$, whose fiber over a smooth curve $X$ is the space of holomorphic one forms ($H^0(X,K_X)$), or equivalently the dual of $H^1(X,\mathcal{O}_X)$.

On the moduli space $\smu{\bz}{0}$ we can define two Hodge-like bundles, according to whether we focus on the base or on the cover curve. The former is however a natural subbundle of the latter, as we shall see in an instant.

\subsubsection{Connected covers}
By the Riemann Hurwitz formula, given an \'{e}tale, connected $\bZ_3$-cover $E\rightarrow C$, the genus of $E$ is $h=3g-2$. There is a natural forgetful morphism  
\be
\smu{\bz}{0}^{conn} \longrightarrow \overline{M}_h,
\ee
and we can define the Hodge bundle on $\smu{\bz}{0}$ by pulling  back $\mathbb{E}^h$ via this morphism. The group action on the covers induces a $\mathbb{Z}_3$ action on $\mathbb{E}^h$, which gives a decomposition 
\be
\mathbb{E}^h= \mathbb{E}_1\oplus \mathbb{E}_\omega \oplus \mathbb{E}_{\bar\omega}
\ee
 into eigenbundles (with respect to the action of the primitive generator of the group). Here $\omega$ and $\bar\omega$ are nontrivial cube roots of unity and denote the corresponding eigenvalues.
 
 The fibers of $\mathbb{E}_1$ are $\bZ_3$-invariant forms, i.e. forms pulled-back from the base curve. It follows that the rank of $\mathbb{E}_1$ is $g$.
 By symmetry arguments the ranks of $\mathbb{E}_\omega$ and $\mathbb{E_{\bar\omega}}$ are $g-1$.
We denote by $\lambda_{i,\omega}$(resp. $\lambda_{i,\bar\omega}$)  the $i$-th Chern class of $\mathbb{E}_\omega$ (resp. $\mathbb{E_{\bar\omega}}$).

\subsubsection{Disconnected covers}
In the case
\be
p: \bigsqcup_{1}^3 C  \longrightarrow C,
\ee
the Hodge bundle corresponding to forms on the cover curves is  a rank $3g$ bundle:  three copies of the Hodge bundle pulled back from $\overline{M}_g$. Keeping track of the $\bZ_3$ action, $\mathbb{E}^{3g}$ is naturally identified with the tensor product of $\mathbb{E}^g$ with the standard representation of $\bZ_3$. The eigenbundles are each a copy of $\mathbb{E}^g$.

\subsection{Invariants of $\orb$ and Hodge integrals}

It might seem deceiving that we discussed at lenght the \GW\ theory of $\bz$ when really we are interested in $\orb$. In fact, typically one cannot even define \GW\ invariants for a non-compact target space, as the moduli space of stable maps  is itself non-compact. When a space $\mathfrak{X}$ admits a torus action with compact fixed locus ${F}$, Bryan and Pandharipande (\cite{bp:bps}) \textit{define} the invariants of $\mathfrak{X}$ via localization: the \GW\ theory of $\mathfrak{X}$ is thus reduced to the \GW\ theory of ${F}$ ``corrected'' by the euler class of an obstruction  (virtual)\footnote{In general it should really be considered as an element in $K$-theory.} bundle constructed from the normal bundle $N_{{F}/\mathfrak{X}}$ (see \cite[Section 2.2]{bp:tlgwtoc}).

A three dimensional torus $(\bC^\ast)^3$ acts naturally on $\bC^3$, and this action descends to the quotient. The only fixed point for the action is the image of the origin, which is a copy of $\bz$. Therefore:
\be
\langle \ \rangle^g_{({\orb},0)} = \int_{\smu{\bz}{0}} e(-R^\bullet\pi_\ast f^\ast(N_{\bz/\orb})).
\ee
\begin{Rem}
Note that the orbifold $\orb$ contains no compact curve classes, therefore the only invariants correspond to constant maps ($\beta=0$).
\end{Rem}
The normal bundle to the image of the origin consists of three copies of a line bundle  denoted $L_\omega$: it descends from a copy of $\bC$ with a non-trivial action of $\bZ_3$. In the world of orbifolds this is an essential feature: the fibers of $R^i\pi_\ast f^\ast (L_\omega)$ over a curve $X$ are not the full $H^i(X,\mathcal{O}_X)$, but only the $\bar\omega$ eigenspace.

Therefore:
\be
R^1\pi_\ast f^\ast(L_\omega)= (\mathbb{E}_\omega)^\vee ,
\ee
and
\be
R^0\pi_\ast f^\ast(L_\omega) =
\left\{
\begin{array}{ll}
\mathcal{O} & \mbox{on\ } \smu{\bz}{0}^{disc} \\
 & \\
0           & \mbox{on\ } \smu{\bz}{0}^{conn}
\end{array}
\right.
\ee
Finally, we are able to express our \GW\ invariants as Hodge integrals:
\begin{eqnarray}
\langle \ \rangle^g_{({\orb},0)} & = &  \frac{1}{t_1t_2t_3}\int_{{\smu{\bz}{0}}^{disc}} e(((\mathbb{E}_\omega)^\vee)^3) \nonumber \\
 &  &+ \int_{{\smu{\bz}{0}}^{conn}}e(((\mathbb{E}_\omega)^\vee)^3)  \nonumber \\
 & = &  \frac{1}{3 t_1t_2t_3}\int_{\overline{M}_g} e((\mathbb{E}^\vee)^3) \label{disc}\\
 &  &+ (-1)^{g-1}\int_{{\smu{\bz}{0}}^{conn}} \lambda_{g-1,\omega}^3 . \label{conn}
\end{eqnarray}

\begin{Rem}
Contribution (\ref{disc}) is a ``classical'' Hodge Integral on the moduli space of stable curves, computed by Faber and Pandharipande in the late '90s (\cite{fp:lsahiittr}, \cite{f:algo}). Contribution (\ref{conn}) is a new and interesting creature, for which we are currently seeking a systematic approach. In low genus one can use ad hoc methods to show that this contribution vanishes.
\end{Rem}

\subsection{Tools for the computation}
The invariants in genus $2$ and $3$ are computed making use of the following classical results.
\begin{description}
		\item[Mumford relation \cite{m:taegotmsoc}] \begin{eqnarray}
c_t(\mathbb{E} \oplus \mathbb{E}^\nu )=  1.
	\label{mumf}\end{eqnarray}
	\item[G-Mumford relation \cite{bgp:crc}] \begin{eqnarray}
c_t(\mathbb{E}_\omega\oplus (\mathbb{E}^\nu)_\omega )= c_t(\mathbb{E}_\omega\oplus (\mathbb{E}_{\bar\omega})^\nu)= 1.
	\label{Gmumf}
	\end{eqnarray}
	\item[Faber-Pandharipande computation \cite{fp:lsahiittr}]
	\begin{eqnarray}
\int_{\overline{M}_g}\lambda_g\lambda_{g-1}\lambda_{g-2}= \frac{1}{2(2g-2)!}\frac{|B_{2g-2}|}{2g-2}\frac{|B_{2g}|}{2g}.
	\label{hodge}
\end{eqnarray}
\end{description}
Here, $c_t$ denotes the Chern polynomial and $B_n$ the $n^{th}$ Bernoulli number.

\section{Invariants of $\orb$ with $g>1$}

In this section we perform some computations of Gromov-Witten invariants of $\orb$ with $g > 1$, using the Hodge integral approach developed in the previous section.

\subsection{$g=2$}

Let us start by computing the genus $2$ unmarked Gromov-Witten invariant of $\orb$.

\subsubsection{Vanishing of (\ref{conn})}
In this case both $\mathbb{E}_\omega$ and $\mathbb{E}_{\bar\omega}$ are line bundles. Integral (\ref{conn}) is:
\be
-\int \lambda_{1,\omega}^3.
\ee

Relation (\ref{Gmumf})yields:
\begin{itemize}
	\item $\lambda_{1,\omega}=\lambda_{1,{\bar\omega}}$;
	\item $\lambda_{1,\omega}\lambda_{1,{\bar\omega}}=0$ .
\end{itemize}
This immediately shows the vanishing of our desired integral.

\subsubsection{Computation of (\ref{disc})}

Integral \eqref{disc} in this case is:
\be
\frac{1}{3t_1t_2t_3}\int_{\overline{M}_2} (\lambda_2 -\lambda_1t_1+t_1^2)(\lambda_2 -\lambda_1t_2+t_2^2)(\lambda_2 -\lambda_1t_3+t_3^2).
\ee
Setting the weights to be Calabi-Yau ($t_1+t_2+t_3=0$), we obtain the following weight independent expression:
\be
\frac{1}{3} \int_{\overline{M}_2} -\lambda_1^3 + 3 \lambda_2\lambda_1 =\frac{1}{3} \int_{\overline{M}_2} \lambda_2\lambda_1,
\ee
where the last equality follows from the application of Mumford's relation (\ref{mumf}) that tells us
that $2\lambda_2=\lambda_1^2$. Using formula (\ref{hodge}), we get:
\be
\begin{imp}\displaystyle{
\langle\ \rangle^2_{(\orb,0)}= \frac{1}{17280}}.
\end{imp}
\ee

\subsection{$g=3$}

We now compute the genus $3$ unmarked Gromov-Witten invariant of $\orb$.

\subsubsection{Vanishing of (\ref{conn})}

In this case the vanishing of (\ref{conn}) is only slightly more elaborate. We want to compute:
\be
A=\int \lambda_{2,\omega}^3 =\int \lambda_{2,\bar\omega}^3 \mbox{\ \ (by symmetry)}.
\ee
Relation (\ref{Gmumf}) gives us:
\begin{description}
	\item[a)]$\lambda_{1,\omega}=\lambda_{1,\bar\omega}=\alpha$,
 \item[b)]$\alpha^2=\lambda_{2,\omega}+\lambda_{2,\bar\omega}$,
 \item[c)]$\alpha\lambda_{2,\omega}=\alpha\lambda_{2,\bar\omega}$,
\item[d)]$\lambda_{2,\omega}\lambda_{2,\bar\omega}=0$.
 \end{description}

Using some elementary algebra and all of the relations above:
\begin{align}
2A&=\int \lambda_{2,\omega}^3+\lambda_{2,\bar\omega}^3=
\int (\lambda_{2,\omega}+\lambda_{2,\bar\omega})(\lambda_{2,\omega}^2+\lambda_{2,\bar\omega}^2) \notag\\
&=\int \alpha^2(\lambda_{2,\omega}^2+\lambda_{2,\bar\omega}^2)= 2\int \alpha^2\lambda_{2,\omega}\lambda_{2,\bar\omega}=0.
\end{align}

\subsubsection{Computation of (\ref{disc})}
The computation  here is identical to genus $2$. With Calabi-Yau weights, and by formula (\ref{hodge}):
\be
\begin{imp}\displaystyle{
\langle \ \rangle^3_{(\orb,0)}= -\frac{1}{3} \int_{\overline{M}_3} \lambda_3\lambda_2\lambda_1=-\frac{1}{4354560}}.
\end{imp}
\ee

\subsection{Higher genus}

Starting with $g=4$, there is no reason why the contribution from the connected covers should vanish. In fact, the prediction from physics, which we will describe in the next section, does not match (\ref{disc}). Invariants with insertions can also be expressed in terms of \zhh, whose structure is still completely unexplored. 
In collaboration with Charles Cadman and Arend Bayer, the second author is attempting a systematic approach of $\bZ_3$ Hodge integrals in higher genus. Currently two avenues are being pursued:
\begin{itemize}
	\item evaluating via localization integrals on auxiliary moduli spaces as a mean to produce relations between \zhh. This approach is similar in spirit to \cite{cc:c3z3}. 
	\item using stacky Grothendieck-Riemann-Roch and the natural covering map between $\smu{\bz}{0}^{conn}$ and $\overline{M}_g$ in order to express \zhh\ in terms of polynomials in tautological classes on $\overline{M}_g$. Such gadgets can then be evaluated through the use of Witten's conjecture, implemented for example in Faber's algorithm \cite{f:algo}.
\end{itemize}

\section{The physics}

In this section we review the calculation of \cite{ABK}.
We first discuss relevant features of the two main ingredients in the calculation, namely mirror symmetry and topological string theory, and then move on to the actual calculation of Gromov-Witten invariants of $\orb$. Good references on mirror symmetry include the two books \cite{CoxKatz, Hori}, while topological string theory is explored in detail in the book \cite{Marino}.

\subsection{Mirror symmetry at large radius}

To start with, we recall the usual local description of mirror symmetry at large radius. The main characters are:
\begin{itemize}
\item $(X,Y)$: a mirror pair of Calabi-Yau threefolds;
\item $\CM(Y)$: a suitable compactification of the complex structure moduli space of $Y$;
\item $\CKM(X)$: a suitable compactification of the complexified K\"ahler moduli space of $X$ --- the so-called \emph{stringy} or \emph{enlarged} K\"ahler moduli space.
\end{itemize}
Mirror symmetry provides a local isomorphism, called the \emph{mirror map}, between $\CKM (X)$ and $\CM (Y)$, which maps a neighborhood of a maximally unipotent boundary point $q_0 \in \CM (Y)$ to a neighborhood of a corresponding large radius point $p_0 \in \CKM (X)$. Moreover, mirror symmetry tells us that the mirror map lifts to an isomorphism between the \emph{A-model} amplitudes at $p_0 \in \CKM (X)$, and the \emph{B-model} amplitudes at $q_0 \in \CM (Y)$.

But what are the A- and B-model amplitudes? Start with a theory --- a non-linear sigma model --- of maps $f : \Sigma \to M$ from Riemann surfaces $\Sigma$ to a Calabi-Yau threefold $M$. There are two ways of twisting this sigma model to obtain topological theories, namely the A- and the B-model. The A-model does not depend on complex moduli, while the B-model is independent of K\"ahler moduli.

\subsubsection{The A-model}
The A-model on $X$ becomes a theory of holomorphic maps $f: \Sigma \to X$, which can be reformulated in terms of Gromov-Witten invariants of the target space $X$. In the neighborhood of $p_0 \in \CKM (X)$, the A-model genus $g$ amplitudes $F_g$ become generating functionals for the unmarked genus $g$ Gromov-Witten invariants $\langle \quad \rangle^g_{(X, \beta)}$ of $X$, that is
\be
\label{GW}
F_g = \sum_{\beta \in H_2 (X)} \langle \quad \rangle^g_{(X, \beta)} Q^{\beta},
\ee
where
\be
Q^{\beta} = \re^{2 \pi i \int_{\beta} \omega},
\ee
and $\omega$ is a complexified K\"ahler class of $X$.

\subsubsection{The B-model}
The B-model on $Y$ localizes on constant maps, and becomes a theory of variations of complex structures of the targe space $Y$. As opposed to their A-model cousins, the B-model amplitudes do not afford a simple mathematical description. Nevertheless, the genus $0$ amplitude can be determined by special geometry,\footnote{See \cite{Freed} for a mathematical exposition of special geometry.} and corresponds to the so-called prepotential. The genus $1$ amplitude can be defined in terms of Ray-Singer torsion of $Y$. For the higher genus amplitudes, one can use the holomorphic anomaly equations of \cite{BCOV} --- which may be understood as some sort of higher genus generalization of special geometry ---  to reconstruct the amplitudes recursively in the neighborhood of $q_0 \in \CM (Y)$, up to an unknown holomorphic function at each genus depending on a finite number of constants. External data, such as boundary conditions, must be used to fix these functions.

Since the A-model on $X$ is mirror to the B-model on $Y$, one can use the B-model point of view to compute the Gromov-Witten theory of the mirror $X$. The two main ingredients entering in the calculation are:
\begin{itemize}
\item the mirror map near the large radius point;
\item a framework to compute the B-model amplitudes near $q_0$, such as special geometry and the holomorphic anomaly equations.
\end{itemize}
This was the strategy used by Candelas et al \cite{Candelas} to compute the number of rational curves in the quintic threefold, which was extended to higher genus in \cite{BCOV}.

\subsection{Global mirror symmetry and orbifold points}

So far we only gave a local description of mirror symmetry, near a large radius point of $\CKM(X)$. However, from a physics point of view, mirror symmetry should be global, in the sense that $\CKM(X)$ should be globally isomorphic to $\CM (Y)$, and similarly for the A- and the B-model amplitudes.

Generically, the stringy K\"ahler moduli space $\CKM(X)$ has a rather complicated structure, which goes beyond the K\"ahler cone of $X$. However, when $X$ is toric, $\CKM (X)$ is also toric and is easily described by the \emph{secondary fan}  associated to $X$ (see for instance \cite{CoxKatz}, section 3.4 and chapter 6, for a more precise discussion). Roughly speaking, $\CKM (X)$ is obtained by gluing along common walls the K\"ahler cones of threefolds birationally equivalent to $X$. Some of these cones correspond to smooth threefolds related to $X$ by flops; each such cone then contains a large radius point, which is mapped by mirror symmetry to a corresponding maximally unipotent boundary point in $\CM (Y)$. However, some other patches correspond to ``non-geometric phases", by which we mean that they are obtained from $X$ by contracting some cycles. In particular, we will be interested in the case where $\CKM (X)$ comprises a patch which contains an \emph{orbifold point} $p_{orb} \in \CKM(X)$, where some cycles of $X$ are contracted to yield an orbifold $\overline X$. This orbifold point is mapped on the mirror side to a point of \emph{finite monodromy} $q_{orb} \in \CM(Y)$, around which monodromy of the periods is finite.\footnote{Here, for simplicity, we implicitly assumed that $\CKM(X)$ and $\CM(Y)$ are one-dimensional, which will be the case for the orbifold $\orb$.}

Our aim is now to study mirror symmetry in the neighborhood of the points $p_{orb}$ and $q_{orb}$. First, one needs to define an \emph{orbifold mirror map}, which identifies these two neighborhoods, and should lift to an isomorphism of the A- and the B-model amplitudes near these points. The relation between A-model amplitudes and Gromov-Witten theory is still valid near $p_{orb}$; namely, the A-model genus $g$ amplitudes now become generating functionals for the genus $g$ orbifold Gromov-Witten invariants of $\overline X$. Hence, our goal is to use the B-model around $q_{orb}$ to compute the orbifold Gromov-Witten invariants of $\overline X$ via the orbifold mirror map. As in the traditional large radius calculation, the essence of the calculation boils down to two ingredients:
\begin{itemize}
\item the orbifold mirror map near the orbifold point;
\item a framework to compute the B-model amplitudes near $q_{orb}$.
\end{itemize}
Let us look at both of these items a little closer.

\subsection{The orbifold mirror map}

\label{mirrormap}

\subsubsection{Large radius point}
At large radius, the mirror map can be described as follows. $H^2 (X, {\mathbb C})$ is spanned by
\be
t_1 T_1 + \ldots + t_r T_r,
\ee
where $T_1, \ldots, T_r$ is a basis of generators for the cone $\sigma$ containing the large radius point $p_0 \in \CKM(X)$ corresponding to $X$. The complexified K\"ahler parameters $t_1, \ldots, t_r$ parameterize $\CKM(X)$ near $p_0$. On the mirror side, as is standard in special geometry we parameterize $\CM(Y)$ using periods of the holomorphic volume form $\Omega$ on $Y$. Choose a symplectic basis of three-cycles $A^I, B_J \in H_3 (Y)$, with $I,J=0, \ldots, r$, and define the periods
\be
\omega^I = \oint_{A^I} \Omega, \qquad {\partial {\CF} \over \partial \omega^{I} } = \oint_{B_I} \Omega,
\ee
where ${\CF}$ is the prepotential. The periods are solutions of the Picard-Fuchs equations, with the following properties. In terms of coordinates $q_i$, $i=1, \ldots, r$ centered at the maximally unipotent boundary point $q_0 \in \CM(Y)$, there is a unique period which is holomorphic, say $\omega^0$, and $r$ periods have logarithmic behavior,
\be
\omega^i = {\omega^0 \over 2 \pi i} \log (q_i) + {\mathcal O}(q), \qquad i = 1, \ldots, r.
\ee
There are $r$ other periods which are quadratic in the logarithm, and one is cubic. The mirror map is then given by 
\be
(t_1, \ldots, t_r) \mapsto {1 \over \omega^0} (\omega^1, \ldots, \omega^r).
\ee
Note that when $X$ and $Y$ are noncompact,\footnote{See for instance \cite{Hosono} for a more precise discussion of special geometry and periods of a noncompact Calabi-Yau threefold $Y$.} the mirror map is simplified by the fact that $\omega^0 = 1$, hence the $t^i$ are directly identified with the logarithmic periods $\omega^i$.

What is important to note here is that the mirror map was fixed by finding:
\begin{enumerate}
\item a canonical basis for the cohomology group $H^2 (X, {\mathbb C})$ at the large radius point $p_0 \in \CKM(X)$;
\item a basis of solutions of the Picard-Fuchs equations (periods) around the maximally unipotent boundary point $q_0 \in \CKM(Y)$ with the required leading behavior.
\end{enumerate}
The second point can also be understood in terms of monodromy properties of the periods. Under monodromy around $q_0$ the logarithmic periods behave as
\be
\omega^i \mapsto \omega^i + 1,
\ee
while on the A-model side the amplitudes are given as an expansion in terms of the exponentiated parameters $Q_i = \re^{2 \pi i t_i}$, see \eqref{GW}. The $Q_i$'s are then invariant under the shift $t_i \mapsto t_i + 1$, which implies that the amplitudes are invariant under monodromy around $q_0$.

\subsubsection{Orbifold point}

To fix the mirror map around the orbifold point $p_{orb} \in \CKM(X)$ we follow the lessons of the previous section. What we want is:
\begin{enumerate}
\item a canonical basis for the orbifold cohomology of $\overline X$;
\item a basis of solutions of the Picard-Fuchs equations near $q_{orb} \in \CM(Y)$ such that the amplitudes are invariant under the finite monodromy around $q_{orb}$.
\end{enumerate}
As simple as it looks, we will see that these two conditions are sufficient to fix unambiguously the orbifold mirror map for simple orbifolds such as $\orb$, up to a scale factor. A prescription equivalent to condition (2) will be to match the representation theoretic data in the orbifold cohomology ring of $\overline X$ to the action of the finite monodromy on the periods.

\subsection{B-model at the orbifold point}

The next item that we need is a formalism to compute the B-model amplitudes near $q_{orb} \in \CM(Y)$. This is provided by the holomorphic anomaly equations of \cite{BCOV}.

Recall that at $q_{orb}$ there is a basis of periods $\omega^i$ which is selected by the orbifold mirror map. As usual the genus $0$ amplitude $F_0$ is simply given by the prepotential $\CF$ of special geometry. For the higher genus amplitudes $F_g$, $g \geq 1$, one can solve the holomorphic anomaly equations near $q_{orb}$ to obtain the following recursive system:
\be
\label{recursive}
F_g = h_g - \Gamma_g \left[ E^{i j}, {\partial \over \partial \omega^{i_1}} \cdots {\partial \over \partial \omega^{i_n}} F_{r < g} \right],
\ee
where $\Gamma_g$ is a functional depending on the derivatives of the lower genus amplitudes $F_{r < g}$ with respect to the periods $\omega^i$, and on the ``propagator" 
\be
E^{i j} = {\partial F_1 \over \partial \tau_{i j}},
\ee
with $\tau_{i j}$ the period matrix:
\be
\tau_{i j} = {\partial^2 {\mathcal F} \over \partial \omega^i \partial \omega^j} = {\partial^2 F_0 \over \partial \omega^i \partial \omega^j}.
\ee 
The $h_g$ are undetermined functions, depending on a finite number of constants. As an example, the genus $2$ functional is given by
\begin{align}
\Gamma_2=&
{E}^{ij}\Big({1\over 2} {\partial_i \partial_j} {F}_1  + {1\over 2} {\partial_i}{F}_1 
{\partial_j} {F}_1\Big) \notag\\
&+   {E}^{ij} E^{kl}\Big( {1\over 2} \partial_i {F}_1  {\partial_j \partial_k \partial_l} {F}_0
+ {1\over 8}  {\partial_i \partial_j \partial_k \partial_l} { F}_0 \Big)\notag\\
&+ {E}^{ij} {E}^{kl} {E}^{mn}\Big({1\over 8} {\partial_i \partial_j \partial_k}
{F}_0 {\partial_l \partial_m \partial_n} {F}_0\notag\\
&\qquad \qquad\qquad\qquad+  {1\over 12} {\partial_i \partial_k \partial_m} {F}_0
{\partial_j \partial_l \partial_n} {F}_0\Big), 
\end{align}
where we used the notation
\be
\partial_i F_k = {\partial F_k\over \partial \omega^i} .
\ee
We refer the reader to \cite{ABK} for the explicit iterative derivation of this recursive system, which is perhaps easier understood in terms of wavefunction properties of the topological string partition function. We note that the holomorphic anomaly equations can also be solved by direct integration using modular properties of the amplitudes, see \cite{GKMW}.

As mentioned earlier, the equations \eqref{recursive} are not complete, in the sense that they cannot be used alone to reconstruct recursively the amplitudes $F_g$, since the holomorphic functions $h_g$ are undetermined. Hence, the system must be supplemented by additional data, such as boundary conditions, to fix the $h_g$'s. 

What kind of additional data can we use at the orbifold point? Well, the simple realization of \cite{ABK} is that we in fact do not need any new data! Indeed, a crucial point is that the $h_g$ are holomorphic functions, which are globally defined all over the moduli space $\CM(Y)$. Hence, if we know the amplitudes at a large radius point $q_0 \in \CM(Y)$, we can fix the $h_g$ and use them, in conjunction with \eqref{recursive}, to compute the amplitudes at the orbifold point $q_{orb} \in \CM(Y)$.

\subsection{Strategy}

\begin{figure}
\begin{center}
\includegraphics[width=4in]{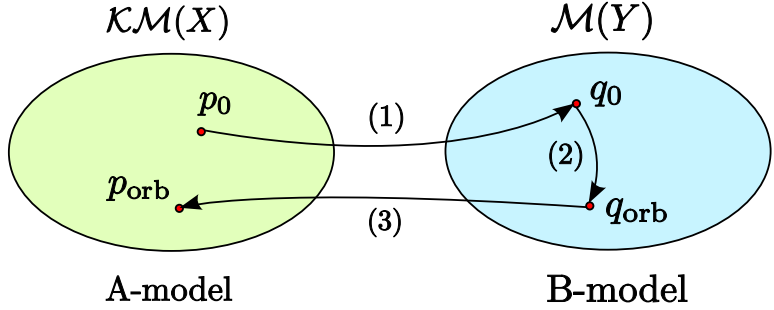}
\end{center}
\caption{A schematic illustration of our strategy to compute orbifold Gromov-Witten invariants.}
\label{f:strategy}
\end{figure}

Our strategy to compute orbifold Gromov-Witten invariants should now be clear. Consider a smooth Calabi-Yau threefold $X$ for which the compactified K\"ahler moduli space $\CKM(X)$ contains an orbifold point $p_{orb} \in \CKM(X)$ corresponding to an orbifold $\overline X$. 
We first determine the mirror maps near the large radius point $p_0 \in \CKM(X)$ and the orbifold point $p_{orb} \in \CKM(X)$, using the principles of section \ref{mirrormap}. The calculation then proceeds in three steps, which are illustrated in figure \ref{f:strategy}.
\begin{enumerate}
\item We compute the generating functionals of Gromov-Witten invariants of $X$, using for instance the topological vertex \cite{AKMV, LLLZ} if $X$ is toric, or localization of Hodge integrals. These are mapped by mirror symmetry at large radius to the B-model amplitudes near $q_0 \in \CM(Y)$.
\item From these amplitudes we fix the holomorphic functions $h_g$, which are valid all over the moduli space and can be used to compute the B-model amplitudes at $q_{orb}$ through the recursion \eqref{recursive}. 
\item Finally, we use the orbifold mirror map to extract the orbifold Gromov-Witten invariants of $\overline X$ from the B-model amplitudes at $q_{orb}$.\footnote{We note here that there is an alternative strategy to compute the orbifold amplitudes, which combines modular --- or wavefunction --- properties of the amplitudes and the symplectic transformation between the periods canonically chosen by the mirror maps at the large radius point and the orbifold point. This was the approach emphasized in \cite{ABK}, where it was shown to be equivalent to the procedure outlined here.}
\end{enumerate}

\section{The physics computation}

We now turn to the calculation of orbifold Gromov-Witten invariants of $\orb$.

\subsection{Mirror symmetry}

The orbifold $\orb$ has a unique crepant resolution, which is the (noncompact) toric Calabi-Yau threefold 
\be
X = {\mathcal O} ( - 3) \to {\mathbb P}^2,
\ee
often called local ${\mathbb P}^2$ in the physics literature. The stringy K\"ahler moduli space $\CKM(X)$ is one-dimensional, and includes two distinct patches; one of which contains the large radius point $p_0$ of $X$, and the other contains an orbifold point $p_{orb}$ where the ${\mathbb P}^2$ is contracted to zero size, yielding the orbifold $\overline X = \orb$.

Following the standard procedure of \cite{HV}, the mirror threefold $Y$ can be described as follows. Let $w, w' \in {\mathbb C}$, and $x,y \in {\mathbb C}^*$. Then $Y$ is the noncompact threefold
\be
Y = \{ w w' = y^2 + y (1 + x) + q x^3 \},
\ee
where $q$ is a coordinate on $\CM(Y)$ centered at the large radius point $q_0 := \{q = 0\} \in \CM(Y)$. That is, $Y$ is a conic fibration over ${\mathbb C}^* \times {\mathbb C}^*$, where the fiber degenerates to two lines over the one-parameter family of Riemann surfaces
\be
\Sigma (q) = \{ y^2 + y (1 + x) + q x^3 = 0 \},
\ee
which has genus $1$ and three punctures. The point of finite monodromy $q_{orb} \in \CM(Y)$ is located at $q \to \infty$. A natural coordinate centered at $q_{orb}$ is
\be
\psi = {(-1)^{1/3} \over 3} q^{-1/3},
\ee
as can be read off from the secondary fan. Note that under ${\mathbb Z}_3$-monodromy around $q_{orb} = \{ \psi = 0 \}$, $\psi$ undergoes
\be
\psi \mapsto \re^{2 \pi i / 3} \psi.
\ee

\subsection{The orbifold mirror map}

The first ingredient that we need to fix is the mirror map near the orbifold point $p_{orb} \in \CKM(X)$. In order to do so, we start by solving the Picard-Fuchs equations near $q_{orb} \in \CM(Y)$. Following the work of Chiang, Klemm, Yau and Zaslow \cite{CKYZ}, we know that the Picard-Fuchs differential operator that annihilates the periods is given by, in terms of the coordinate $\psi$ centered at $q_{orb}$:
\be
{\mathcal D}_\psi= \psi^3 \Theta_\psi^3 - ( \Theta_\psi-2)(\Theta_\psi -1)\Theta_\psi,
\ee
with $\Theta_\psi = \psi \partial_\psi$. ${\mathcal D}_\psi \Pi^{orb}=0$ can be solved with  techniques from \cite{gkz:gkz}; a solution vector is given by $\Pi^{orb} = (1, B_1(\psi), B_2(\psi) )$ with
\be
B_k(\psi)= \sum_{n \geq 0} {(-1)^{3 n+k+1} \psi^{3 n + k} \over (3 n + k)!} \left({ \Gamma \left( n + {k \over 3} \right) \over \Gamma \left( {k \over 3} \right) } \right)^3.
\ee

As described in section \ref{mirrormap}, to get the orbifold mirror map we need to find linear combinations of the solutions above that are mapped to a basis for the orbifold cohomology of $\orb$. The orbifold cohomology $H^*_{\rm orb} (\orb)$ has basis ${\bf 1}_0$, ${\bf 1}_{1/3}$ and ${\bf 1}_{2/3}$, where the ${\bf 1}_{r/3}$ are components of the inertia stack corresponding to the elements $[r]$ of ${\mathbb Z}_3$. The basis elements have degrees 
\be
\deg({\bf 1}_0) = 0, \qquad \deg( {\bf 1}_{1/3}) = 2, \qquad \deg ( {\bf 1}_{2/3} ) = 4.
\ee 
Hence $H^*_{\rm orb}(\orb)$ is spanned by
\be
\sigma_0 {\bf 1}_0 + \sigma_1 {\bf 1}_{1/3} +\sigma_2 {\bf 1}_{2/3}.
\ee
The orbifold mirror map will be given by mapping $\sigma_1$ to an appropriate combination of $1$, $B_1(\psi)$ and $B_2(\psi)$.

Recall that monodromy around $q_{orb}$ is given by $\psi \mapsto \re^{2 \pi i /3} \psi$, which implies
\be
(1, B_1(\psi), B_2(\psi) ) \mapsto (1, \re^{2 \pi i /3} B_1(\psi), \re^{4 \pi i /3} B_2(\psi) ).
\ee
But ${\bf 1}_{1/3}$ corresponds to the element $[1] \in {\mathbb Z}_3$, or, in terms of third roots of unity, to $\re^{2 \pi i /3}$. Thus, it is clear that $\sigma_1$ must be mapped to $B_1(\psi)$ directly, up to an overall scale factor. More precisely, we claim that the mirror map is given by
\be
\label{mirrormapc3}
(\sigma_1, \sigma_2) = (B_1(\psi), B_2(\psi) ).
\ee
Another way of arguing for this mirror map is by computing the genus $0$ amplitude, as we do next. Up to scale, the above mirror map is the only map that yields a genus $0$ amplitude which is invariant under orbifold monodromy. Note that this is also the mirror map that was proved in \cite{CCIT}.

\subsection{Genus $0$ amplitude}

Before computing the genus $0$ amplitude, let us clarify the relation between the A-model amplitudes and Gromov-Witten theory at the orbifold point. At large radius, the genus $g$ A-model amplitudes become generating functionals for genus $g$ Gromov-Witten invariants $\langle \quad \rangle_{(X, \beta)}^g$ in homology classes $\beta \in H_2 (X, {\mathbb Z})$, with no insertions. At the orbifold point, $\orb$ contains no compact curve, hence the only invariants correspond to constant maps $\beta = 0$. However, the A-model amplitudes now become generating functionals for orbifold Gromov-Witten invariants with marked points, more precisely
\be
F_g^{orb} = \sum_{n=0}^{\infty}{1 \over n!} \langle ({\bf 1}_{1/3})^n \rangle^g_{(\orb, 0)} \sigma_1^{n}.
\ee 
Note that the unmarked ($n=0$) invariants are only well-defined for $g \geq 2$. Moreover, only contributions with $n \in 3 {\mathbb Z}$ are non-zero, which ensures that the amplitudes are invariant under orbifold monodromy.

To compute the genus $0$ amplitude, we use the fact that it is given by the prepotential of special geometry, which is defined by\footnote{The unusual factor of $-3$ here comes from the fact that since $Y$ is noncompact, it is not possible to find a symplectic basis of three-cycles; instead, the A- and the B-cycles have intersection number $-3$.}
\be
\sigma_2 = -3 {\partial \CF^{orb} \over \partial \sigma_1 }.
\ee
${\mathcal F}^{orb}$ gives the genus $0$ orbifold Gromov-Witten potential $F_0^{orb}$ of $\orb$. Integrating $\sigma_2$, we get:
\be
F_0^{orb} = \sum_{k=1}^\infty {1 \over (3k)!} \langle ({\bf 1}_{1/3})^{3 k} \rangle^{g=0}_{(\orb, 0)} \sigma_1^{3 k}
\ee
with the invariants $N_{0,k}:= \langle ({\bf 1}_{1/3})^{3 k} \rangle^{0}_{(\orb, 0)}$:
\be
N_{0,1} = {1\over 3}, \;\; N_{0,2} = -{1\over 3^3}, \;\; N_{0,3} = {1\over 3^2},
 \;\; N_{0,4}=-{1093 \over 3^6}, \ldots
\ee
 Agreement with the mathematical computation  of the genus $0$ amplitude fixes the normalization of the mirror map \eqref{mirrormapc3}. 

\subsection{Higher genus amplitudes}

To extract the higher genus amplitudes of $\orb$, we need to compute the holomorphic functions $h_g$ at each genus $g$. This can be done easily at large radius, by first computing the A-model amplitudes through the topological vertex, and then mapping them to the B-model side using the usual mirror map at large radius. 
We obtain, for the marked invariants:
\be
F_g^{orb} = \sum_{k=1}^\infty {1 \over (3k )!} \langle ({\bf 1}_{1/3})^{3 k} \rangle^{g}_{(\orb, 0)}\sigma^{3k}
\ee
with the numbers $N_{g,k} := \langle ({\bf 1}_{1/3})^{3 k} \rangle^{g}_{(\orb, 0)}$:
\begin{center}
\begin{small}
\begin{tabular}{|c|rrrr|}\hline
$g$
&$k=1$
&$2$
&$3$
&$4$\\
\hline
$0$
&${1 \over 3}$
&$-{1 \over 3^3}$
&${1 \over 3^2}$
&$-{1093 \over 3^6}$\\[0.2cm]
$1$
&$0$
&${1\over 3^5}$
&$-{14\over 3^5}$
&${13007\over 3^8}$\\[0.2cm]
$2$
&${1\over 2^4\cdot 3^4 \cdot 5}$
&$-{13\over 2^4\cdot3^6}$
&${20693\over 2^4\cdot 3^8 \cdot 5}$
&$-{12803923\over 2^4 \cdot 3^{10}\cdot 5}$\\[0.2cm]
$3$
&$-{31\over 2^5 3^5 5 \cdot 7}$
&${11569\over 2^5 3^9 5\cdot 7}$
&$-{2429003 \over 2^5 3^{10} 5\cdot 7}$
&${871749323\over 2^4 3^{11}5\cdot 7}$\\[0.2cm]
$4$
&${313\over 2^7 3^9 5^2}$
&$-{1889\over 2^7 3^{9} }$
&${115647179\over 2^6 3^{13} 5^2 }$
&$-{29321809247\over 2^8 3^{12} 5^{2}}$\\[0.2cm]
$5$
&$-{519961\over 2^9 3^{11} 5^2 7 \cdot 11}$
&${196898123\over 2^9 3^{12} 5^2 7\cdot 11}$
&$-{339157983781\over 2^9 3^{14} 5^2 7\cdot 11}$
&${78658947782147 \over 2^9 3^{16} 5 \cdot 7}$\\[0.2cm]
$6$
&${ 14609730607\over 2^{12} 3^{13} 5^3 7^2  11}$
&$-{258703053013\over 2^{10} 3^{15} 5^1 7^2 11}$
&${2453678654644313\over 2^{12} 3^{14} 5^3 7^2 11}$
&$-{40015774193969601803 \over 2^{11} 3^{18} 5^3 7^2 11}$\\[0.2cm]
\hline
\end{tabular}
\end{small}
\end{center}

The unmarked invariants ($n=k=0$) for $g \geq 2$ (these are not well-defined for $g=0,1$) can also be calculated, and read
\begin{gather}
N_{2,0} = {-1 \over 2160} + { \chi(X) \over 5760}, \;\; N_{3,0} = {1\over 544320} - {\chi(X) \over 1451520}, \\ N_{4,0} = -{7\over 41990400}+{\chi(X) \over 87091200},
N_{5,0} = {3161 \over 77598259200} - {\chi(X) \over 2554675200}, \ldots\notag 
\end{gather}
where $\chi(X)$ is the ``Euler number" of $X = {\mathcal O}(-3) \to {\mathbb P}^2$. 

A little more should be said about the unmarked invariants. To compute these invariants, we first needed the degree $0$ unmarked invariants at large radius, that is the invariants $\langle \quad \rangle^g_{(X,0)}$ for constant maps to $X = {\mathcal O}(-3) \to {\mathbb P}^2$, which give the second term in each of the expressions above. These invariants were computed by Faber and Pandharipande: 
\be
\langle \quad \rangle^g_{(X,0)} = (-1)^g \chi(X) \int_{\overline{M}_g} \lambda_g \lambda_{g-1} \lambda_{g-2}, \qquad g \geq 2,
\ee
where we use the notation of section 1. The integral here is precisely the Hodge integral \eqref{hodge}.
Even though talking about the Euler characteristic of a noncompact threefold might make some differential geometers cringe, we observe that any vector bundle retracts to its zero section. Therefore, $\chi(X)=\chi(\Proj)=3$, and we obtain:

\be
N_{2,0} = {1 \over 17280}, \qquad N_{3,0} = - {1 \over 4354560},
\ee
which match perfectly the results obtained earlier via Hodge integrals. 

To end this section, let us mention that although the calculation of the unmarked invariants here is relatively similar to the Hodge integral calculation performed earlier (in particular the use of Faber-Pandharipande's formula), it is fundamentally different. Indeed, as noticed in section 2, for $g \geq 4$ the second integral \eqref{conn} should not vanish anymore, and the direct Hodge integral calculation necessitates an understanding of these new ${\mathbb Z}_3$-Hodge integrals. However, for the physics calculation, only the Faber-Pandharipande standard Hodge integral is needed, since the corrections come from the functions $h_g$ and the recursive formula \eqref{recursive}.

\newcommand{\etalchar}[1]{$^{#1}$}
\def\cprime{$'$}


\begin{thebibliography}{CDLOGP91}

\bibitem[ABK06]{ABK}
M. Aganagic, V. Bouchard and A. Klemm.
\newblock Topological strings and (almost) modular forms.
\newblock arXiv:hep-th/0607100, 2006.

\bibitem[ACV03]{Abramovich_Corti_Vistoli}
D.~Abramovich, A.~Corti and A.~Vistoli.
\newblock Twisted bundles and admissible covers.
\newblock {\em Comm. Algebra}, 31(8):3547--3618, 2003  [arXiv:math.AG/0106211].
\newblock Special issue in honor of Steven L. Kleiman.

\bibitem[AGV06]{Abramovich_Graber_Vistoli}
D.~Abramovich, T.~Graber and A.~Vistoli.
\newblock Gromov-{W}itten theory of {D}eligne-{M}umford stacks.
\newblock arXiv:math.AG/0603151, 2006.

\bibitem[AKMV05]{AKMV}
M. Aganagic, A. Klemm, M. Mari\~no and C. Vafa.
\newblock The topological vertex.
\newblock {\em Commun. Math. Phys.}, 254:425--478, 2005 [arXiv:hep-th/0305132].

\bibitem[BC07]{bc:c3z3}
A. Bayer and C. Cadman.
\newblock Quantum cohomology of $[\mathbb{C}^n/\mu_r]$.
\newblock arXiv:0705.2160 [math.AG], 2007.

\bibitem[BCOV94]{BCOV}
M.~Bershadsky, S.~Cecotti, H.~Ooguri and C.~Vafa.
\newblock Kodaira-Spencer theory of gravity and exact results for quantum
  string amplitudes.
\newblock {\em Commun. Math. Phys.}, 165:311--428, 1994 [arXiv:hep-th/9309140].

\bibitem[BG06]{bg:crc}
J. Bryan and T. Graber.
\newblock The crepant resolution conjecture.
\newblock arXiv:math.AG/0610129, 2006.

\bibitem[BG07a]{bg:crchhi}
J. Bryan and A. Gohlampour.
\newblock Root Systems and the Quantum Cohomology of ADE resolutions
\newblock arXiv:0707.1337 [math.AG], 2007.

\bibitem[BG07b]{bg:root}
J. Bryan and A. Gohlampour.
\newblock Hurwitz-Hodge integrals, the $E_6$ and $D_4$ root systems, and the crepant
  resolution conjecture.
\newblock arXiv:0708.4244 [math.AG], 2007.

\bibitem[BGP05]{bgp:crc}
J. Bryan, T. Graber and R. Pandharipande.
\newblock The orbifold quantum cohomology of $\mathbb{C}^2/\mathbb{Z}_3$ and
  {Hurwitz-Hodge} integrals.
\newblock arXiv:math.AG/0510335, 2005.

\bibitem[BP01]{bp:bps}
J. Bryan and R. Pandharipande.
\newblock {BPS} states of curves in {Calabi-Yau} $3$-folds.
\newblock Geometry and Topology, 5:287--318, 2001 [arXiv:math.AG/0009025].

\bibitem[BP04]{bp:tlgwtoc}
J. Bryan and R. Pandharipande.
\newblock The local {Gromov-Witten} theory of curves.
\newblock arXiv:math.AG/0411037, 2004.

\bibitem[CC07]{cc:c3z3}
C. Cadman and R. Cavalieri.
\newblock Gerby localization, $\mathbb{Z}_3$-{H}odge integrals and the {GW}
  theory of $\orb$.
\newblock arXiv:0705.2158 [math.AG], 2007.

\bibitem[CCIT06]{ccit:wc}
T. Coates, A. Corti, H. Iritani and H.-H. Tseng.
\newblock Wall-crossings in toric Gromov-Witten theory I: Crepant examples.
\newblock arXiv:math.AG/0611550, 2006.

\bibitem[CCIT07a]{CCIT}
T. Coates, A. Corti, H. Iritani and H.-H. Tseng.
\newblock Computing genus-zero twisted Gromov-Witten invariants.
\newblock arXiv:math.AG/0702234, 2007.

\bibitem[CCIT07b]{ccit:crcas}
T. Coates, A. Corti, H. Iritani and H.-H. Tseng.
\newblock The crepant resolution conjecture for type {A} surface singularities.
\newblock arXiv:0704.2034 [math.AG], 2007.

\bibitem[CDLOGP91]{Candelas}
P. Candelas, X. de~la~Ossa, P.~S. Green and L. Parkes.
\newblock A pair of Calabi-Yau manifolds as an exactly soluble superconformal
  theory.
\newblock {\em Nucl. Phys.}, B359:21--74, 1991.

\bibitem[CK99]{CoxKatz}
D.~A. Cox and S. Katz.
\newblock {\em Mirror Symmetry and Algebraic Geoemetry}, vol.~68 of {\em
  Mathematical Surveys and Monographs}.
\newblock AMS, 1999.

\bibitem[CKYZ99]{CKYZ}
T.~M. Chiang, A.~Klemm, S.-T. Yau and E.~Zaslow.
\newblock Local mirror symmetry: Calculations and interpretations.
\newblock {\em Adv. Theor. Math. Phys.}, 3:495--565, 1999 [arXiv:hep-th/9903053].

\bibitem[CR02]{cr:ogwt}
W. Chen and Y. Ruan.
\newblock Orbifold {G}romov-{W}itten theory.
\newblock In {\em Orbifolds in mathematics and physics (Madison, WI, 2001)},
  volume 310 of {\em Contemp. Math.}, pages 25--85. Amer. Math. Soc.,
  Providence, RI, 2002 [arXiv:math.AG/0103156].

\bibitem[CR04]{cr:nctoo}
W. Chen and Y. Ruan.
\newblock A new cohomology theory of orbifolds.
\newblock {\em Comm. Math. Phys.}, 248(1):1--31, 2004 [arXiv:math.AG/0004129].

\bibitem[DLOFS02]{OFS}
  X.~de la Ossa, B.~Florea and H.~Skarke.
\newblock D-branes on noncompact Calabi-Yau manifolds: K-theory and monodromy.
\newblock  {\em Nucl.\ Phys.\  B}, 644:170, 2002 [arXiv:hep-th/0104254].

\bibitem[DG00]{DG}
  D.~E.~Diaconescu and J.~Gomis.
\newblock  Fractional branes and boundary states in orbifold theories.
\newblock  {\em JHEP} 0010:001, 2000  [arXiv:hep-th/9906242].

\bibitem[Fab99]{f:algo}
C. Faber.
\newblock Algorithms for computing intersection numbers on moduli spaces of
  curves, with an application to the class of the locus of jacobians.
\newblock {\em New trends in algebraic geometry (Warwick 1996),London Math.
  Soc. Lecture Note Ser.}, 264:93--109, 1999 [arXiv:alg-geom/9706006].

\bibitem[FP00]{fp:lsahiittr}
C.~Faber and R.~Pandharipande.
\newblock Logarithmic series and {H}odge integrals in the tautological ring.
\newblock {\em Michigan Math. J.}, 48:215--252, 2000 [arXiv:math.AG/0002112].
\newblock With an appendix by Don Zagier, Dedicated to William Fulton on the
  occasion of his 60th birthday.

\bibitem[Fre99]{Freed}
  D.~S.~Freed.
\newblock  Special Kaehler manifolds.
\newblock {\em Commun.\ Math.\ Phys.},  203:31, (1999) [arXiv:hep-th/9712042].

\bibitem[Gil07]{g:crc}
W. D. Gillam.
\newblock The crepant resolution conjecture for 3-dimensional flags modulo an
  involution.
\newblock arXiv:0708.0842 [math.AG], 2007.

\bibitem[GKZ94]{gkz:gkz}
I.~M. Gel{\cprime}fand, M.~M. Kapranov, and A.~V. Zelevinsky.
\newblock {\em Discriminants, resultants, and multidimensional determinants}.
\newblock Mathematics: Theory \& Applications. Birkh\"auser Boston Inc.,
  Boston, MA, 1994.

\bibitem[GKMW07]{GKMW}
 T. W. Grimm, A. Klemm, M. Mari\~no and M. Weiss.
\newblock Direct integration of the topological string. 
\newblock arXiv:hep-th/0702187, 2007.

\bibitem[HKK{\etalchar{+}}03]{Hori}
K. Hori, S. Katz, A. Klemm, R. Pandharipande, R.
  Thomas, C. Vafa, R. Vakil and E. Zaslow.
\newblock {\em Mirror Symmetry}, volume~1 of {\em Clay Mathematics Monographs}.
\newblock AMS, 2003.

\bibitem[HV00]{HV}
K. Hori and C. Vafa.
\newblock Mirror symmetry.
\newblock arXiv:hep-th/0002222, 2000.

\bibitem[Hos04]{Hosono}
  S.~Hosono.
\newblock Central charges, symplectic forms, and hypergeometric series in local mirror symmetry.
\newblock  arXiv:hep-th/0404043, 2004.

\bibitem[LLLZ04]{LLLZ}
J. Li, C.-C.~M. Liu, K. Liu, and J. Zhou.
\newblock A mathematical theory of the topological vertex.
\newblock arXiv:math.AG/0408426, 2004.

\bibitem[Mar05]{Marino}
M.~Mari\~no.
\newblock {\it Chern-Simons theory, matrix models, and topological strings.}
\newblock Oxford, UK: Clarendon, 2005.

\bibitem[Mum83]{m:taegotmsoc}
D. Mumford.
\newblock Toward an enumerative geometry of the moduli space of curves.
\newblock {\em Arithmetic and Geometry}, II(36):271--326, 1983.

\bibitem[Rua01]{r:crc}
Y. Ruan.
\newblock Cohomology ring of crepant resolutions of orbifolds.
\newblock arXiv:math.AG/0108195, 2001.

\end{thebibliography}
\end{document}